\documentclass[12pt]{article}
\usepackage{bbm}
\usepackage{epsfig}
\usepackage{amssymb,amsmath,amscd}
\usepackage{graphicx,cite}
\usepackage{subfigure, graphicx}

\newcommand{\EOP} { \hfill $\Box$ }

\newtheorem{theorem}{Theorem}[section]
\newtheorem{prop}[theorem]{Proposition}

\begin{document}
\title{\bf  Computation of Maxwell's equations on Manifold using DEC}

\author{
Zheng Xie$^1$\thanks{E-mail: lenozhengxie@yahoo.com.cn }~~~~Yujie
Ma$^2$\thanks{ E-mail: yjma@mmrc.iss.ac.cn \  This work is partially
supported by  CPSFFP (No. 20090460102), NKBRPC (No. 2004CB318000),
and NNSFC (No. 10871170) }
\\{\small
$1.$ Center of Mathematical Sciences, Zhejiang University
(310027),China}
\\ {\small $2.$ Key Laboratory of Mathematics Mechanization,}
\\ {\small  Chinese Academy of Sciences,  (100090), China}}

\date{}

\maketitle

\begin{abstract}
In this paper, the method of discrete exterior calculus for
numerically solving Maxwell's equations in space manifold and the
time is discussed, which is  a kind of lattice gauge theory. The
analysis of its   stable condition and error is also accomplished.
This algorithm has been implemented on C++ plateform for simulating
TE/M waves in vacuum.

\end{abstract}


\vskip 0.2cm \noindent {\bf PACS(2010):} 41.20.Jb, 02.30.Jr,
02.40.Sf, 02.60.Cb.

\section{Introduction}

Computational electromagnetism is concerned with the numerical study
of Maxwell's equations. The Yee scheme is known as finite difference
time domain  and is one of the most successful numerical methods,
particularly in the area of microwave problems\cite{yee}. It
preserves important structural features of Maxwell's
equations\cite{bondeson,clemens,gross,hairer,nicolet}. Bossavit et
al present the
 Yee-like scheme and
 extend   Yee scheme to unstructured grids. This scheme  combines the
 best attributes of the finite element method (unstructured grids)
  and Yee scheme (preserving geometric structure)\cite{bossavit1,bossavit2}.
 Stern et al \cite{stern} generalize the Yee scheme to
unstructured grids not just in space, but  in 4-dimensional
spacetime by discrete exterior calculus(DEC)\cite{
whitney,arnold,dimakis,desbrun,meyer,leok,hyman,hiptmair,wise,novikov,auchmann}.
This relaxes the need to take uniform time steps.

In this paper,  we  generalize the Yee  scheme to the
 discrete space manifold and the time. The spacetime manifold used here
 is split as a product of $1$D time and $2$D or $3$D space manifold.
 The  space manifold can be approximated by triangular  and tetrahedrons depending on dimension,
  and the time by segments. So the spacetime manifold is approximated by prism lattice,
 on which the discrete Lorentz metric can be defined.
\begin{itemize}
  \item[1.]
    With the technique of discrete exterior calculus, the $\mathbb{R}$ value
  discrete connection,
  curvature and  Bianchi identity are defined on prim lattice.   With discrete
variation of an inner product of   discrete  $1-$forms and their
dual forms, the discrete source equation and continuity equation are
derived.

  \item[2.] Those equations
compose the discrete Maxwell's equations in vacuum case, which  just
need the local information of triangulated manifold such as length
and area.
\end{itemize}The discrete Maxwell's equations here can be re-grouped
  into two sets  of explicit iterative schemes for
  TE and TM waves, respectively.
Those schemes can directly use acute triangular, rectangular,
regular polygon and their combination, which has been implemented on
C++ plateform to simulate the electromagnetic waves propagation and
interference on manifold.

\section{DEC for Maxwell's equations}

Maxwell's equations can be simply expressed once the language of
exterior differential forms is used. The electric and magnetic
fields are jointly described by a curvature $2-$form $F$ in a 4-D
spacetime manifold. The Maxwell's equations  reduce to the Bianchi
identity and the source equation
$$\begin{array}{lll}dF=0 && d\ast F=\ast J\end{array}\eqno{(1)}$$
where $d$ denotes the exterior differential operator, $\ast$ denotes
the Hodge star operator, and  1-form $J$ is called the  electric
current form satisfying the continuity equation
$$d\ast J=0.$$

As the exterior derivative is defined on any manifold, the
differential form version of the Bianchi identity makes sense for
any 3D or 4D spacetime manifold, whereas the source equation is
defined if the manifold is oriented and has a Lorentz metric. Now,
we introduce the discrete counterpart of those differential
geometric objects to derive the numerical computational schemes for
Maxwell's equations.

\subsection*{Discrete Lorentz metric }
 The spacetime manifold used here
 is split as a product of $1$D time and $2$D or $3$D space manifold.
 The $2$D or $3$D   space manifold can be approximated by triangular  or tetrahedrons,
  and the time by segments.  The length of edge and area of triangular and volume of tetrahedrons
gives the discrete Riemann metric on space grids. The metric on time
grid is the minus of length  square. The spacetime manifold is
approximated by prism lattice,
 on which the discrete Lorentz metric can be defined as
the product of   discrete  metric on space and time.

\subsection*{Discrete exterior calculus}

 A discrete differential $k$-form, $k \in \mathbb{Z}$, is the
evaluation   of the differential $k$-form on all $k$-simplices. Dual
forms, i.e., forms that we evaluate on the dual cell.  Suppose each
simplex contains its circumcenter. The circumcentric dual cell
$D(\sigma_0)$ of simplex $\sigma_0$ is
 $$ D(\sigma_0):=\bigcup_{\sigma_0\in \sigma_1\in\cdots \in\sigma_r}
 \mathrm{Int}(c(\sigma_0)c(\sigma_1)\cdots c(\sigma_r) ),$$
where   $\sigma_i$ is all the simplices which contains
$\sigma_0$,..., $\sigma_{i-1}$, and $c(\sigma_i)$ is the
circumcenter of $\sigma_i$.

The   two operators in Eqs.(1) can be discretized as follows:
\begin{itemize}
  \item[1.] Discrete exterior differential operator $d$, this operator is the transpose of the
incidence matrix of $k$-cells on $k+1$-cells.
  \item[2.] Discrete Hodge Star $\ast$, the operator     scales the cells by the volumes of the corresponding dual and
primal cells.
\end{itemize}

\subsection*{Discrete connection and  curvature}
Discrete connection $1-$form or gauge field $A$ assigns to each
element
 in the  set of edges $E$ an element of the gauge group $\mathbb{R}$:
$$A: E\rightarrow \mathbb{R}.$$
Discrete curvature $2-$form is the discrete exterior derivative of
the   discrete  connection $1-$form
$$F=dA: P\rightarrow \mathbb{R}.$$ The value of $F$ on
each element in the set of triangular $P$ is   the coefficient of
Holonomy group of this face. The $2-$form  $F$ automatically
satisfies the discrete Bianchi identity
$$d F=0.\eqno{(2)}$$

Note that since the gauge group $\mathbb{R}$ used here  is abelian,
we need not pick a starting vertex for the loop. We may traverse the
edges in any order, so long as  taking orientations into account.

%
%
\subsection*{Discrete Maxwell's equations  }
%

For source case, we need discrete current $1-$form  $J$. Let
$A=\sum\limits_EA_i$  and the Lagrangian functional be
 $$\begin{array}{lll}
L(A,J)&=& -\frac{1}{2}\langle dA,  dA\rangle +\langle A,J\rangle,
\end{array}
$$
where
$$\begin{array}{lll}\langle dA,  dA\rangle &:=& (A)_{1\times |E|}(d)_{|E|\times|F|}
(\ast)_{|F|\times|F|}(d)^{T}_{|F|\times|E|}(A)^{T}_{|E|\times 1}\\
\langle A,J\rangle&:=&(A)_{1\times |E|} (\ast)_{|E|\times|E|}
(J^{T})_{|E|\times 1}
\end{array}
$$

Supposing that there is a variation of $A_i$, vanishing on the
boundary, we have
$$\begin{array}{lll}
\partial_{A_i}L(A,J)&=&\partial_{A_i}(-\frac{1}{2}\langle dA,  dA\rangle +\langle A,J\rangle)\\
&=&\partial_{A_i}(-\frac{1}{2} (A)_{1\times
|E|}(d)_{|E|\times|F|}(*)_{|F|\times|F|}(d)^{T}_{|F|\times|E|}(A)^{T}_{|E|\times
1} \\&&+(A)_{1\times |E|} (*)_{|E|\times|E|} (J^{T})_{|E|\times 1})\\
 &=&   -\frac{1}{2} (0,...,\underbrace{1}_i,...,0)_{1\times
|E|}(d)_{|E|\times|F|}(*)_{|F|\times|F|}(d)^{T}_{|F|\times|E|}(A)^{T}_{|E|\times
1}\\
&&-\frac{1}{2} (A)_{1\times
|E|}(d)_{|E|\times|F|}(*)_{|F|\times|F|}(d)^{T}_{|F|\times|E|}(0,...,\underbrace{1}_i,...,0)^{T}_{|E|\times
1}\\
&&
+(0,...,\underbrace{1}_i,...,0)_{1\times |E|} (*)_{|E|\times|E|} (J^{T})_{|E|\times 1}    \\
 &=& -(0,...,\underbrace{1}_i,...,0)_{1\times
|E|}(d)_{|E|\times|F|}(*)_{|F|\times|F|}(d)^{T}_{|F|\times|E|}(A)^{T}_{|E|\times
1}\\
&&+(0,...,\underbrace{1}_i,...,0)_{1\times |E|} (*)_{|E|\times|E|}
(J^{T})_{|E|\times 1}
\end{array}
$$

The Hamilton's principle of stationary action states that this
variation must equal zero for any such vary of $A_i$, implying the
Euler-Lagrange equations
$$   - (d)_{|E|\times|F|}(*)_{|F|\times|F|}(d)^{T}_{|F|\times|E|}(A)^{T}_{|E|\times
1} +  (*)_{|E|\times|E|} (J^{T})_{|E|\times 1}=0, $$ which is the
discrete  source equation
$$\delta F=J, \eqno{(3)}$$ where
$\delta=\ast^{-1}d^{T}\ast$. Since $(d^{T})^2=0$, the discrete
continuity equation can express as:
$$d^{T}\ast J=0. \eqno{(4)}$$
The equations of discrete Bianchi identity (2),   source equation
(3), and continuity equation (4) are called discrete Maxwell's
equations.

\subsection*{Discrete Gauge transformations}
Discrete gauge transformations are maps
$$A\rightarrow  A+df$$ for any $0-$form or scalar function $f$ on
vertex.
 Since  the discrete exterior derivative maps
$$F\rightarrow  F +d^2f = F,$$
the discrete Maxwell's equations(2-4) are invariant under discrete
gauge transformations.   Since the discrete continuity equation (4)
ensures
$$ \langle df,J\rangle=(f)_{1\times |V|}(d^T)_{|V|\times |E|}(\ast)_{|E|\times|E|}
(J^{T})_{|E|\times 1}=0,
$$ so we have  $$\begin{array}{lll}
L(A+df,J)&=& -\frac{1}{2}\langle d(A+df),  d(A+df)\rangle +\langle
A+df,J\rangle= L(A,J).
\end{array}
$$
That is to say the Lagrangian function   is also invariant under
discrete gauge transformations.

\section{Explicit schemes}

\subsection*{Schemes for TE wave}

The  discrete current $1-$form,  discrete curvature $2-$form and its
dual can be written as
$$J=(-\rho_e dt,
J_e)~~~~F^{n+\frac{1}{2}}=E^{n+\frac{1}{2}}\wedge d {t}+B^{n}
~~~~\ast F^{n}=H^{n } \wedge d {t}-D^{n-\frac{1}{2}},$$
 where  $n$ and $n+\frac{1}{2}$ denote the coordinate of the time,
   $E=\sum\limits_{E}{E}_ie^i$ (electric field) is discrete $1-$form on space,
    $B=\sum
\limits_{P} {B}_i P^i$ (magnetic field) is discrete $2-$form on
space,  $H=\sum\limits_{P} {H}_i
*P^i$ (magnetizing field) is the dual of $B$   on space,
 $D=\sum\limits_{E}  {D}_i
*e^i$(electric displacement field) is the  dual of $E$  on
space,  $\rho_e dt$ (charge density) is the  discrete $1-$form on
time,  $J_e=\sum\limits_{E}J_{ei}e^i$ (electric current density) is
the  discrete $1-$form on space. The discrete Maxwell's equations
can be rewritten as
$$\begin{array}{lll}
d_sB^{n}&=&0\\
d_s E^{n+\frac{1}{2}}\wedge d {t} &=&-d_{ {t}}B^n\\
d^{T}_sD^{n-\frac{1}{2}}&=&\ast (\rho_e dt)^{n-\frac{1}{2}}\\
d^{T}_s H^{n}\wedge d{t} &=&d^{T}_{{t}}D^{n-\frac{1}{2}}+\ast
 {J}^{n}_e , \end{array}$$
where $d_s$, $d^T_{{s}}$ are the restriction of $d$, $d^T$ on space,
and
$$\begin{array}{lll}d_{ {t}} B^n : = \dfrac{B^{n+1}-B^{n}}{\Delta t} \wedge{d
{t}}&&
   d^{T}_{
{t}}D^{n-\frac{1}{2}} := \dfrac{D^{n+\frac{1}{2}}-D^{n-\frac{1}{2}}}
{\Delta t}\wedge d {t}.\end{array} \eqno{(5)}$$
 If  allowing for the
possibility of magnetic charges and current discrete $3-$form
$$\bar{J} =(\rho_m,-{J}_m\wedge d {t}),$$ the symmetric scheme can  be
written as
 $$\begin{array}{lll}
d_sB^{n}&=&   \rho^{n }_m\\
d_s E^{n+\frac{1}{2}}\wedge d {t} &=&-d_{ {t}}B^n-
{J}^{n+\frac{1}{2}}_m\wedge d {t}  \\
d^{T}_sD^{n-\frac{1}{2}}&=&\ast (\rho_e dt)^{n-\frac{1}{2}}\\
d^{T}_s H^{n}\wedge d{t} &=&d^{T}_{{t}}D^{n-\frac{1}{2}}+\ast
 {J}^{n}_e , \end{array}\eqno{(6)}$$
where\begin{itemize}
       \item[~] $ \rho_m=\sum\limits_{Tet}\rho_{mi}T^i$(magnetic charges) is discrete 3-form  on  space,
       \item[~] ${J}_m=\sum\limits_{P}J_{mi}P^i$
 (current) is discrete $2-$form on space.
     \end{itemize}
The compact form of Eqs.(6) can be written as
$$\begin{array}{lll}dF=\bar{J}&& d^T\ast F=\ast J, \end{array}$$
with discrete continuity equations or integrability conditions
$$\begin{array}{lll}
d\bar{J} =0&& d^T\ast J=0. \end{array}$$


\begin{prop}\label{TH:fundmat} If the initial condition satisfies the first and third equations in
Eqs.(6),  the solution of the second and fourth equations in Eqs.(6)
automatically satisfy Eqs.(6).
\end{prop}
\noindent Proof.  Because the dimension of spacetime is $3+1$,
therefore
$$d^{T}_s\ast(\rho_e dt)=0 ~~~~d^{T}_t\ast J_e=0~~~~d_s
\rho_m=0~~~~d_t  J_m\wedge dt=0,$$ and the continuity equations can
be reduced to $$d^T_t(\ast \rho_e dt^{n-\frac{1}{2}})-d^T_s\ast
J^{n}_e=0~~~~~-d_sJ^{n+\frac{1}{2}}_m\wedge dt+d_t\rho^n_m=0.$$ So
we have
$$\begin{array}{lll}d^{T}_td^{T}_sD^{n-\frac{1}{2}}-d^{T}_t\ast
(\rho_e dt)^{n-\frac{1}{2}}&=&-d^{T}_t\ast (\rho_e
dt)^{n-\frac{1}{2}}-d^{T}_s(d^{T}_s H^{n}\wedge
d{t}-\ast J^{n}_e)\\
&=&0
\\d_td_sB^{n}-d_t\rho^{n }_m &=& -d_t\rho^{n }_m+d_s(d_s E^{n+\frac{1}{2}}\wedge d
{t}+ {J}^{n+\frac{1}{2}}_m\wedge d {t}) \\
&=&0.
\end{array}
$$
\EOP

Now we show the scheme (5)  on  the product of 2D discrete space
manifold and time. The second and fourth equations in Eqs.(5) based
on Fig.1 are
\begin{equation*}
\left.
    \begin{aligned}
&\frac{ {D}^{n+\frac{1}{2}}_1 -  {D}^{n-\frac{1}{2}}_1 }{\Delta
t}+{J^{n}_e}_1
=\frac{ {H}^{n}_{1} - {H}^{n}_{2} }{|*e_1|}&\\
&-\frac{ {B}^{n+1}_1 - {B}^{n}_1 } {\Delta t}=\frac {
{E}^{n+\frac{1}{2}}_1 |e_1|+ {E}^{n+\frac{1}{2}}_2 |e_2|+
{E}^{n+\frac{1}{2}}_3 |e_3|}{|P_1|}.&
\end{aligned}
  \right\} \eqno{(7)}
\end{equation*}
 where $|~|$ denotes the measure of forms and dual.
The summation on the right is orient, that is to say, inverse the
orientation of $e_i$, then multiply $-1$ with $\bar{E}_i$. Eqs.(7)
can be implemented on 2D discrete manifold directly(see Fig.1).
Eq.(7) on rectangular gird is just the Yee scheme.
$$
\begin{minipage}{0.99\textwidth}
\begin{center}\includegraphics[scale=0.4]{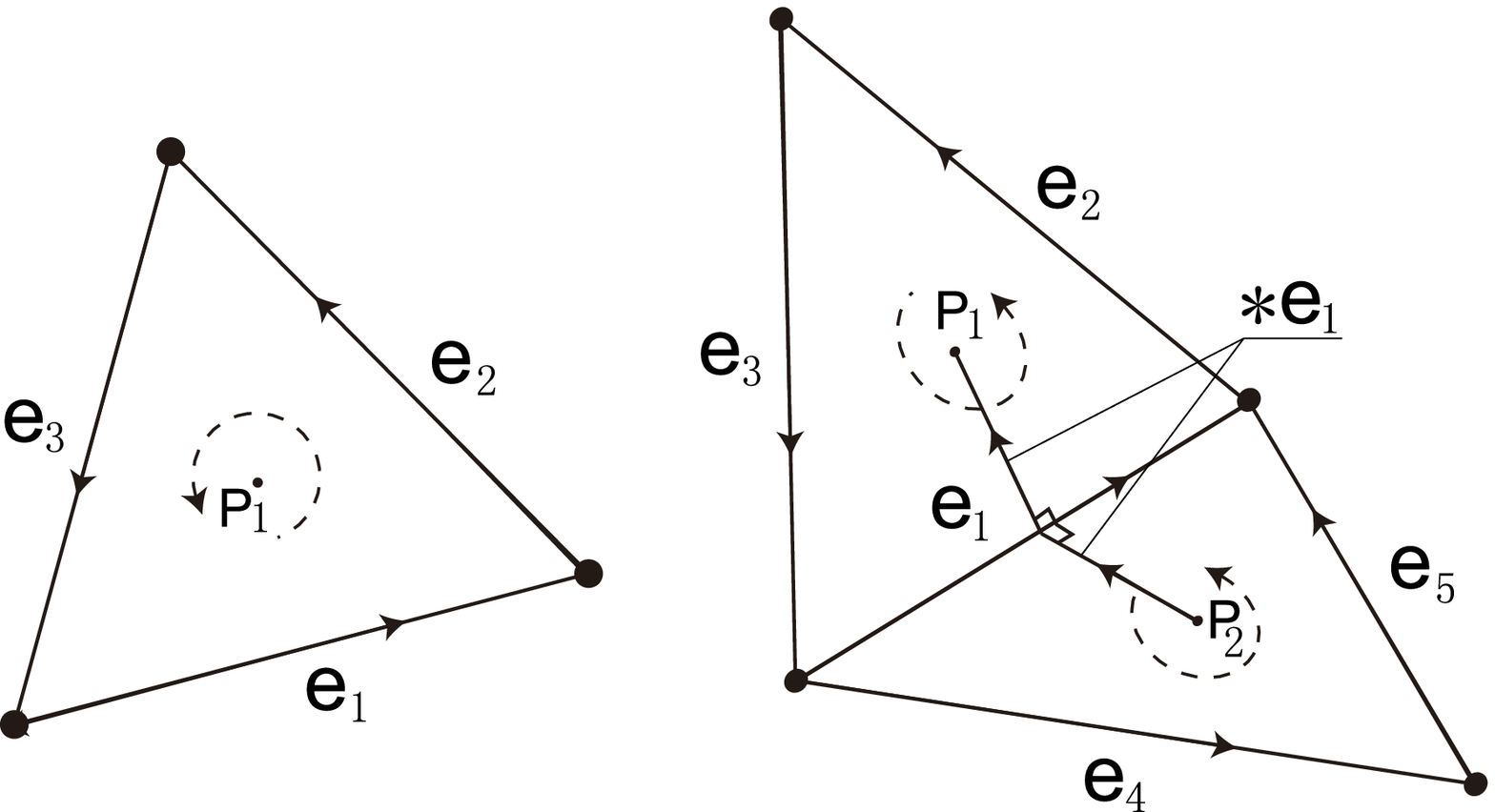}
\end{center}
\begin{center}
{Figure 1: edge and face with  direction  }
\end{center}
\end{minipage}
$$

In the absence of magnetic or dielectric materials, the relations
are simple: $$\begin{array}{lll}D_i=\varepsilon_0 E_i &&
B_i=\mu_0H_i,\end{array}\eqno{(8)}$$ where $\varepsilon_0$ and
$\mu_0$ are two universal constants, called the permittivity of free
space and permeability of free space, respectively. With relations
(8), Eqs.(7) can be rewritten
  into an explicit iterative scheme  for
  TE  wave.
\begin{equation*}
\left.
    \begin{aligned}
&\epsilon_0\dfrac{ {E}^{n+\frac{1}{2}}_1 - {E}^{n-\frac{1}{2}}_1 }
{\Delta t}
+J^{n}_{e1}=\dfrac{ {H}^{n}_{1} - {H}^{n}_{2} }{|*e_1|}&\\
&\mu_0\dfrac{ {H}^{n+1}_1-  {H}^{n}_1 }{\Delta t} =-\dfrac{
{E}^{n+\frac{1}{2}}_1 |e_1|+ {E}^{n+\frac{1}{2}}_2 |e_2|+
{E}^{n+\frac{1}{2}}_3 |e_3|}{|P_1|}&
\end{aligned}
  \right\}\mathrm{TE}\eqno{(9)}
\end{equation*}
The  symmetric TE wave scheme induced from Eqs.(6) can  be written
as follows.
\begin{equation*}
\left.
    \begin{aligned}
&\epsilon_0\dfrac{ {E}^{n+\frac{1}{2}}_1 - {E}^{n-\frac{1}{2}}_1 }
{\Delta t}
+J^{n}_{e1}=\dfrac{ {H}^{n}_{1} - {H}^{n}_{2} }{|*e_1|}&\\
&\mu_0\dfrac{ {H}^{n+1}_1-  {H}^{n}_1 }{\Delta
t}+J^{n+\frac{1}{2}}_{m1} =-\dfrac{ {E}^{n+\frac{1}{2}}_1 |e_1|+
{E}^{n+\frac{1}{2}}_2 |e_2|+ {E}^{n+\frac{1}{2}}_3 |e_3|}{|P_1|}&
\end{aligned}
  \right\}\mathrm{TE}\eqno{(10) }
\end{equation*}

\subsection*{Schemes for TM wave}

If writing
$$\begin{array}{lll}
F^{n+\frac{1}{2}}= H^{n+\frac{1}{2}}\wedge d {t}-D^{n}  &&\ast
F^{n}=-E^{n } \wedge d {t}-B^{n-\frac{1}{2}  }\\
 \bar{J}=(-\rho_e,{J}_e\wedge dt)&& {J}=(-\rho_m dt,{J}_m),
\end{array}$$
 where  $H=\sum\limits_{E} {H}_i
e^i$  is the discrete   $1-$form   on space,  $D=\sum\limits_{P}
{H}_i P^i$  is the discrete   $2-$form  on space,
$E=\sum\limits_{*P} {D}_i
*P^i$  is the dual of $D$  on space,
   $B=\sum\limits_{*E} {B}_i
*e^i$   is the dual of $H$  on space,
 $\rho_e=\sum\limits_{Tet}\rho_{ei}T^i$  is the discrete   $3-$form
on space,  ${J}_e=\sum\limits_{P}J_{ei}P^i$  is the discrete
$2-$form on space,  $\rho_m dt$  is the discrete  $1-$form   on
time,
 ${J}_m=\sum\limits_{E}J_{mi}E^i$  is the discrete   $1-$form   on
space,  the discrete Maxwell's equations can  be rewritten as
$$\begin{array}{lll}
d_sD^{n}&=& \rho^n_e\\
d_s H^{n+\frac{1}{2}}\wedge d {t} &=& d_{ {t}}D^n+J^{n+\frac{1}{2}}_e\wedge dt\\
d^{T}_sB^{n-\frac{1}{2}}&=&\ast (\rho_m  dt)^{n-\frac{1}{2}}\\
d^{T}_s E^{n}\wedge d{t} &=&-d^{T}_{{t}}B^{n-\frac{1}{2}}-\ast
{J}^{n}_m  .
\end{array}\eqno{(11)}$$
\begin{prop}\label{TH:fundmat} If the initial condition satisfies the first and third equations in
Eqs.(11),  the solution of the second and fourth equations in
Eqs.(11) automatically satisfy Eqs.(11).
\end{prop}
\noindent Proof.Because the dimension of spacetime is $3+1$ or
$2+1$, therefore $$d^{T}_s\ast(\rho_m dt)=0~~~d^{T}_t\ast
J_m=0~~~d_s \rho_e=0~~~d_t  J_e\wedge dt=0,$$ and the continuity
equations can be reduced to$$d^{T}_t\ast( \rho_m
dt)^{n-\frac{1}{2}}-d^{T}_s(\ast
J^{n}_m)=0~~~~~d_sJ^{n+\frac{1}{2}}_e\wedge dt-d_t\rho^n_e=0.$$ So
we have
$$\begin{array}{lll}d_td_sD^{n}-d_t
\rho_e^{n}&=&-d_t \rho_e^{n}-d_s(d_s H^{n+\frac{1}{2}}\wedge
d{t}- J^{n+\frac{1}{2}}_e\wedge dt) \\
&=&0
\\d^T_td^T_sB^{n-\frac{1}{2}}-d^T_t\ast (\rho_m  dt)^{n-\frac{1}{2}}  &=& -d^T_t\ast (\rho_m  dt)^{n-\frac{1}{2}}+d^T_s(d^T_s E^{n+\frac{1}{2}}\wedge d
{t}+ \ast {J}^{n}_m )\\
&=&0.
\end{array}
$$
\EOP

 Now we show the scheme (11)  on  the product of 2D discrete space
manifold and time.  The second and fourth equations in Eqs.(11)
based on Fig.1 are
\begin{equation*}
\left.
    \begin{aligned}
&\frac{ {B}^{n+\frac{1}{2}}_1 -  {B}^{n-\frac{1}{2}}_1 }{\Delta
t}+{J^{n}_m}_1
=-\frac{ {E}^{n}_{1} - {E}^{n}_{2} }{|*e_1|}&\\
& \frac{ {D}^{n+1}_1 - {D}^{n}_1 } {\Delta
t}+{J^{n+\frac{1}{2}}_{e1}}=\frac { {H}^{n+\frac{1}{2}}_1 |e_1|+
{H}^{n+\frac{1}{2}}_2 |e_2|+ {H}^{n+\frac{1}{2}}_3 |e_3|}{|P_1|}.&
\end{aligned}
  \right\}\eqno{(12)}
\end{equation*}
 With relations (8),
Eqs.(12) can be rewritten
  into an explicit iterative scheme  for
  TM  wave.
\begin{equation*}
\left.
 \begin{aligned}
    &\epsilon_0\dfrac{ {E}^{n+1}_1- {E}^n_1 }
{\Delta t}+ {J}^{n+\frac{1}{2}}_{e1} =\dfrac{
{H}^{n+\frac{1}{2}}_1|e_1|+ {H}^{n+\frac{1}{2}}_2
|e_2|+ {H}^{n+\frac{1}{2}}_3|e_3|}{|P_1|}&\\
&\mu_0\dfrac{ {H}^{n+\frac{1}{2}}_1 -  {H}^{n-\frac{1}{2}}_1}{\Delta
t}+  {J}^{n}_{m1} =-\dfrac{ {E}^{n}_{1} - {E}^{n}_{2}}{|*e_1|}
&\end{aligned}
  \right\}\mathrm{TM}\eqno{(13)}
\end{equation*}

\subsection*{General schemes }

For real world materials, the constitutive relations are not simple
proportionalities, except approximately. The relations can usually
still be written:
$$\begin{array}{lll}D=\varepsilon  E &&
B=\mu H,\end{array}\eqno{}$$ but $\varepsilon$ and $\mu$ are not, in
general, simple constants, but rather functions. With Ohm's law
$$E=\dfrac{1}{\sigma}J,~~~~~~J_m=\dfrac{1}{\sigma_m}H ,$$where $\sigma$ is the electrical conductivity and $\sigma_m$ is
magnetic conductivity. The DEC schemes can be written as
\begin{equation*}
\left.
    \begin{aligned}
&\epsilon\dfrac{ {E}^{n+\frac{1}{2}}_1 - {E}^{n-\frac{1}{2}}_1 }
{\Delta t}+{\sigma}\dfrac{ {E}^{n+\frac{1}{2}}_1 +
{E}^{n-\frac{1}{2}}_1 }{2}
=\dfrac{ {H}^{n}_{1} - {H}^{n}_{2} }{|*e_1|}&\\
&\mu\dfrac{ {H}^{n+1}_1-  {H}^{n}_1 }{\Delta t}+\sigma_m\dfrac{
{H}^{n+1}_1 + {H}^{n}_1} {2} =-\dfrac{ {E}^{n+\frac{1}{2}}_1 |e_1|+
{E}^{n+\frac{1}{2}}_2 |e_2|+ {E}^{n+\frac{1}{2}}_3 |e_3|}{|P_1|},&
\end{aligned}
  \right\}\mathrm{TE}
\end{equation*}
\begin{equation*}
\left.
 \begin{aligned}
    &\epsilon\dfrac{ {E}^{n+1}_1- {E}^{n }_1 }
{\Delta t}+{\sigma}\dfrac{ {E}^{n+1}_1 + {E}^{n }_1 }{2} =\dfrac{
{H}^{n+\frac{1}{2}}_1
|e_1|+ {H}^{n+\frac{1}{2}}_2 |e_2|+ {H}^{n+\frac{1}{2}}_3 |e_3|}{|P_1|}~~~~~~&\\
&\mu\dfrac{ {H}^{n+\frac{1}{2}}_1 -  {H}^{n-\frac{1}{2}}_1 }{\Delta
t}+\sigma_m\dfrac{ {H}^{n+\frac{1}{2}}_1 + {H}^{n-\frac{1}{2}}_1
}{2} =-\dfrac{ {E}^n_{1} - {E}^n_{2} }{|*e_1|}.&
\end{aligned}
  \right\}\mathrm{TM}
\end{equation*}

\section{Stability, convergence  and accuracy}

\subsection*{Stability}
The Courant-Friedrichs-Lewy condition   is a necessary condition for
convergence while solving certain partial differential equations
numerically.    Now, we find this condition for scheme (10).
Condition for scheme  (13) can be induced in the same way. First, we
decompose DEC algorithm into temporal and spacial eigenvalue
problems.

The temporal eigenvalue problem:
$$\dfrac{\partial^2H^n_0}{\partial t^2}=\Lambda H^n_0$$
It can approximated by difference equation
$$\dfrac{H^{n+1}_0-2H^{n}_0+H^{n-1}_0}{(\Delta t)^2}=\Lambda H^n_0.\eqno{(14)}$$
Supposing $$H^{n+1}_0= H^n_0 \cos( n_1\Delta t)~~~~H^{n-1}_0=H^{n
}_0 \cos( n_2\Delta t)$$  and $$H^{n+1}_0= H^n_0 \sin( n_1\Delta
t)~~~~H^{n-1}_0=H^{n }_0 \sin( n_2\Delta t),$$ and substituting
those into Eq.(14), we obtain
$$\dfrac{\cos( n_1\Delta t)+\cos( n_2\Delta t)-2}{ (\Delta t)^2}  =\Lambda,$$
$$\dfrac{\sin( n_1\Delta t)+\sin( n_2\Delta t)-2}{ (\Delta t)^2}  =\Lambda,$$
therefore
$$-\dfrac{ 4}{ ( \Delta t)^2}  \leq \Lambda\leq 0.$$
This is the   stabile condition for the temporal eigenvalue problem.

The spacial eigenvalue problem:
$$c^2\Delta H=\Lambda H$$
It can be approximated by difference equation (15) based on Fig.2.
  $$
\begin{array}{lll}
   \dfrac{P_{123 }}{c^2}\Lambda H_0 &=& \dfrac{l_{23} }{l_{A0}}(H_A
-H_0 )+\dfrac{l_{34} }{l_{B0} }(H_B -H_0 )+\dfrac{l_{45} }{l_{C0}
}(H_C -H_0 )
\end{array}\eqno{(15)}
$$
$$
\begin{minipage}{0.99\textwidth}
\begin{center}\includegraphics[scale=0.3]{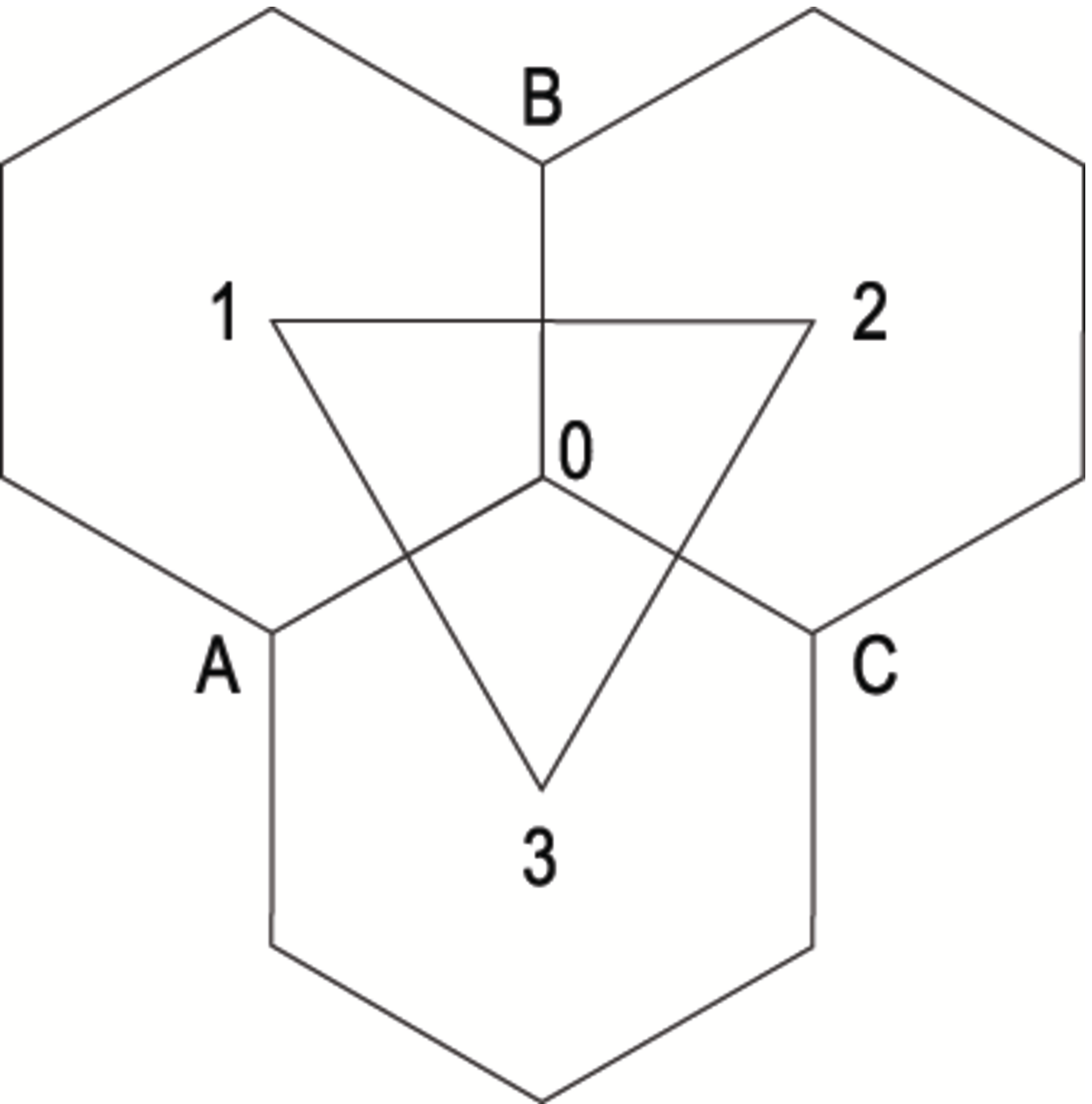}
\end{center}
\begin{center}
{Figure 2: Face and dual face }
\end{center}
\end{minipage}
$$
Let $$H_i= H_0 \cos(c l_{0i})~~~\mathrm{or}~~~H_i=H_0 \sin(c
l_{0i}),$$ and substitute into Eq.(15) to obtain
  $$
  \dfrac{P_{123}}{c^2}\Lambda = \dfrac{l_{23} }{l_{A0}}(\cos(c l_{0A})
-1 )+\dfrac{l_{34} }{l_{B0} }(\cos(c l_{0B}) -1 )+\dfrac{l_{45}
}{l_{C0} }(\cos(c l_{0C}) -1 )
$$
$$
  \dfrac{P_{123}}{c^2}\Lambda =\dfrac{l_{23} }{l_{A0}}(\sin(c l_{0A})
-1 )+\dfrac{l_{34} }{l_{B0} }(\sin(c l_{0B}) -1 )+\dfrac{l_{45}
}{l_{C0} }(\sin(c l_{0C}) -1 ).
$$
So we have
$$-\dfrac{2c^2}{P_{123}}\left(\dfrac{l_{23} }{l_{A0} } +\dfrac{l_{34}
}{l_{B0} } +\dfrac{l_{45} }{l_{C0} }   \right) \leq\Lambda\leq 0.$$
In order to keep the stability of scheme (13), we need
$$-\dfrac{ 2}{ ( \Delta t)^2}\leq -\dfrac{c^2}{P_{123}}\left(\dfrac{l_{23} }{l_{A0} } +\dfrac{l_{34}
}{l_{B0} } +\dfrac{l_{45} }{l_{C0} } \right),\eqno{(16)}$$ or
$$ \Delta t \leq {\mathrm{Mim}}_{P_{123}\in P} \left(\dfrac{1}{c}\sqrt{\dfrac{2P_{123}}{\left(\dfrac{l_{23} }{l_{A0} } +\dfrac{l_{34}
}{l_{B0} } +\dfrac{l_{45} }{l_{C0} }\right)}}\right).$$

\subsection*{Convergence}

By the definition of truncation error, the exact solution
  of Maxwell's equations satisfy the same relation as  DEC scheme
  except for an additional term $O((\Delta t)^2+\Delta t|\ast e|)$.
   This expresses the consistency, and so  convergence for  DEC scheme by  Lax equivalence theorem
   (consistency $+$ stability $=$ convergence).

\subsection*{Accuracy}

The derivative of Maxwell's equations is approximated by first order
difference in schemes (10) and (13). Equivalently, $H$ and $E$ are
approximated by linear interpolation functions. Consulting the
  definition about accuracy of finite volume method, we
can also say that schemes (10) and (13) have first order spacial and
temporal accuracy, and have second order spacial and temporal
accuracy on rectangular grid with equivalent space and time steps.

\section{Implementation}

The DEC algorithm of Maxwell's equations was implemented  in C++
platform. The Fig.3 shows the flowchart of  DEC schemes for
Maxwell's equations.
$$
\begin{minipage}{0.99\textwidth}
\begin{center}\includegraphics[scale=0.7]{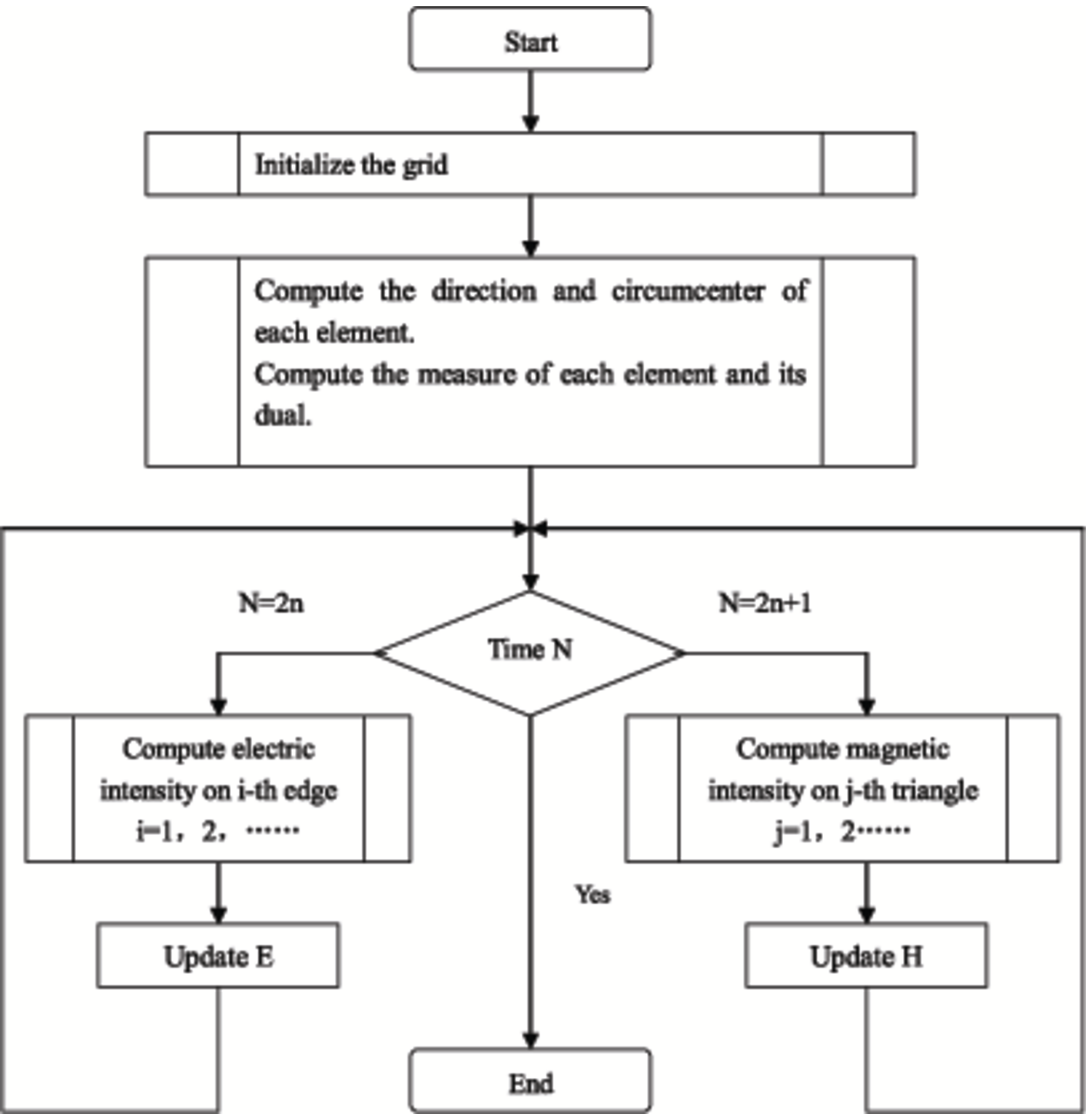}
 \end{center}
\begin{center}
{Figure 3 }
\end{center}
\end{minipage}
$$

In the common practice, not every simulation step needs to be
visualized, especially when the time step size is too small.
  Fig.4
exhibit Gaussian pluses' waveforms simulated by DEC.
$$
\begin{minipage}{0.99\textwidth}
\begin{center}\includegraphics[scale=0.335]{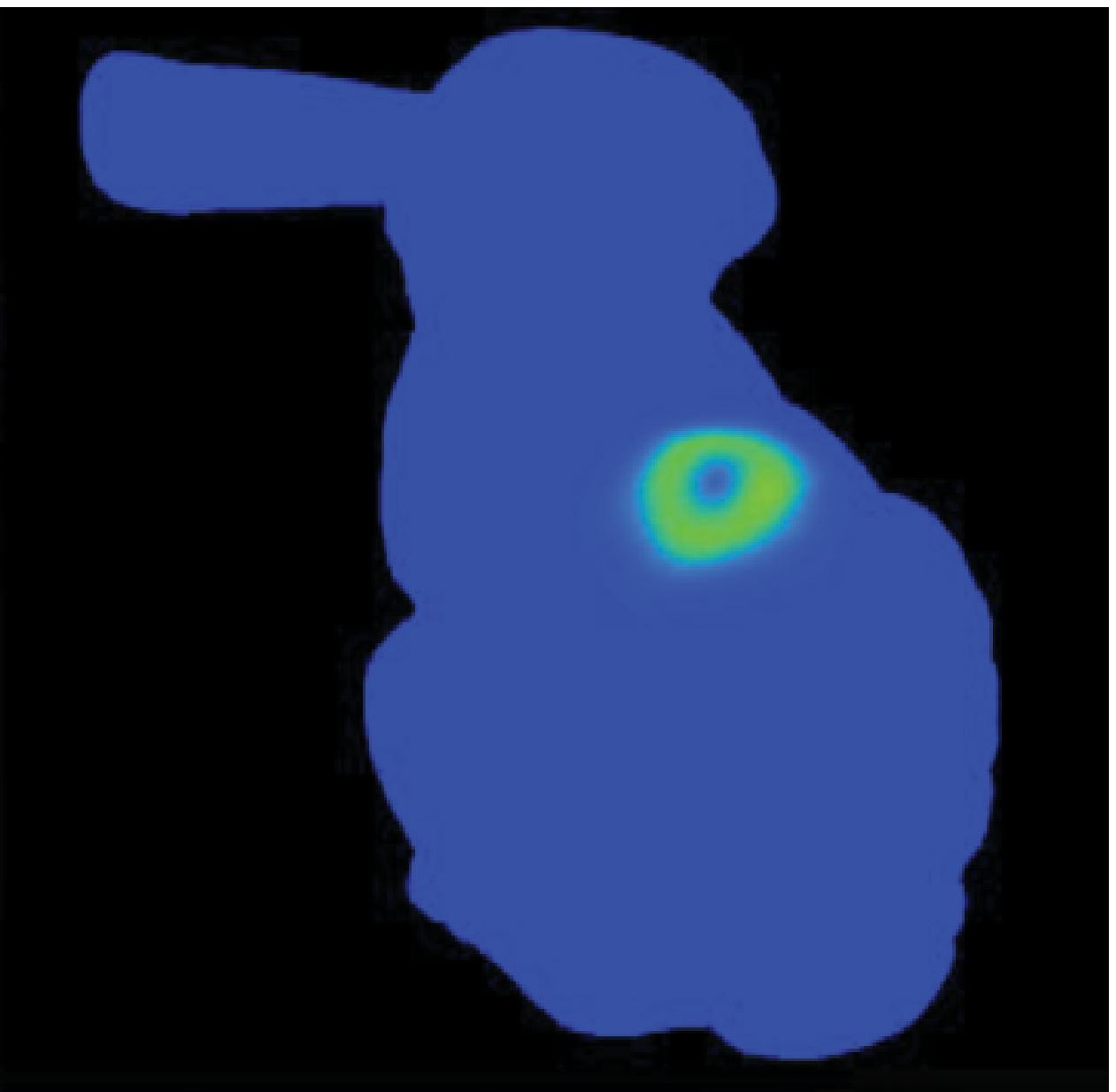}
\includegraphics[scale=0.335]{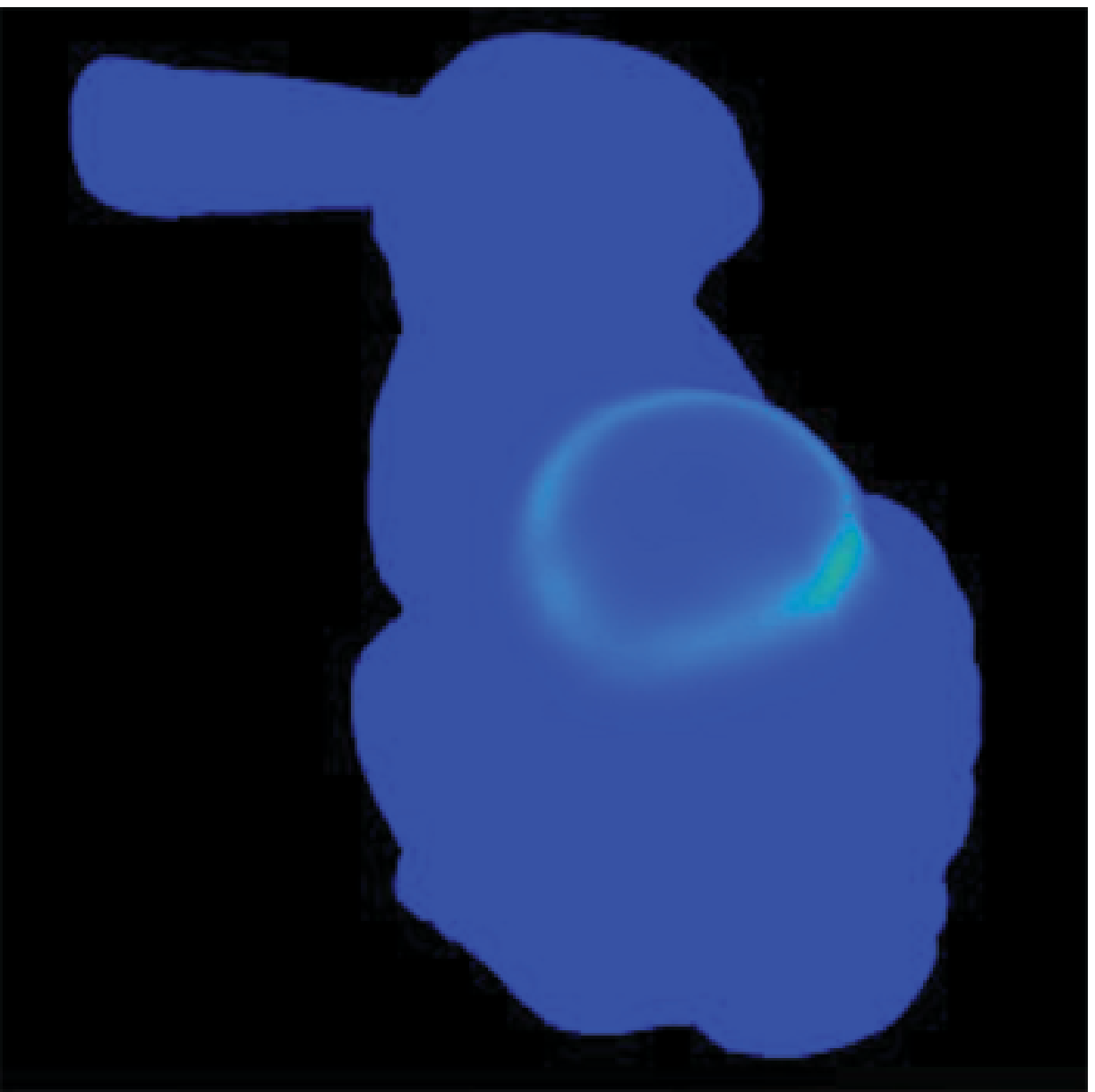}
\includegraphics[scale=0.335]{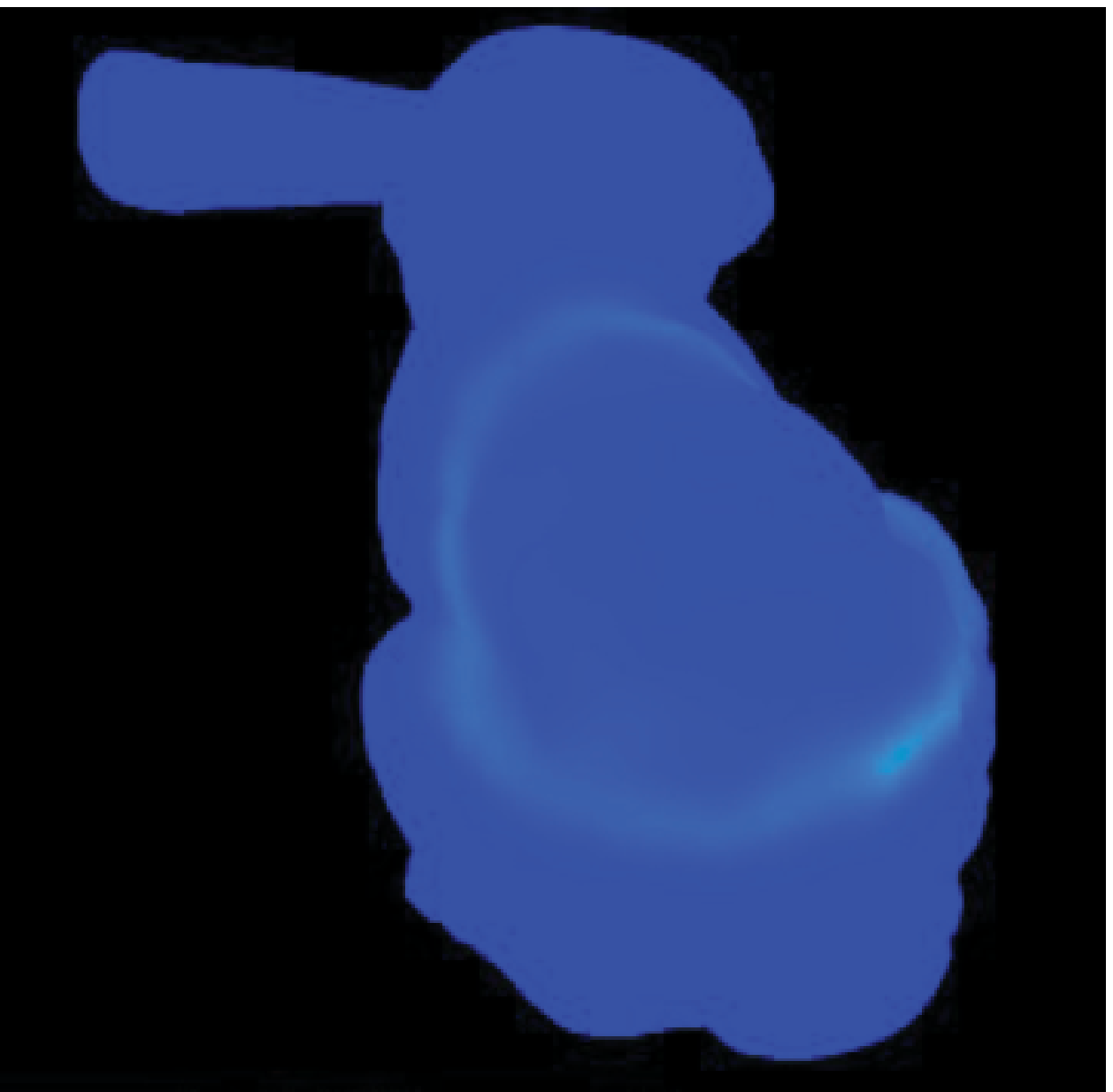}
\end{center}
\begin{center}
{Figure 4:  Simulation of  Gaussian pluse on Stanford bunny by DEC}
\end{center}
\end{minipage}
$$
Fig.5 exhibits two sources Gaussian pluses interference simulated by
DEC. Our algorithm can simulate more complex situation on surface
and  3D-space manifold.
$$
\begin{minipage}{0.99\textwidth}
\begin{center}\includegraphics[scale=0.335]{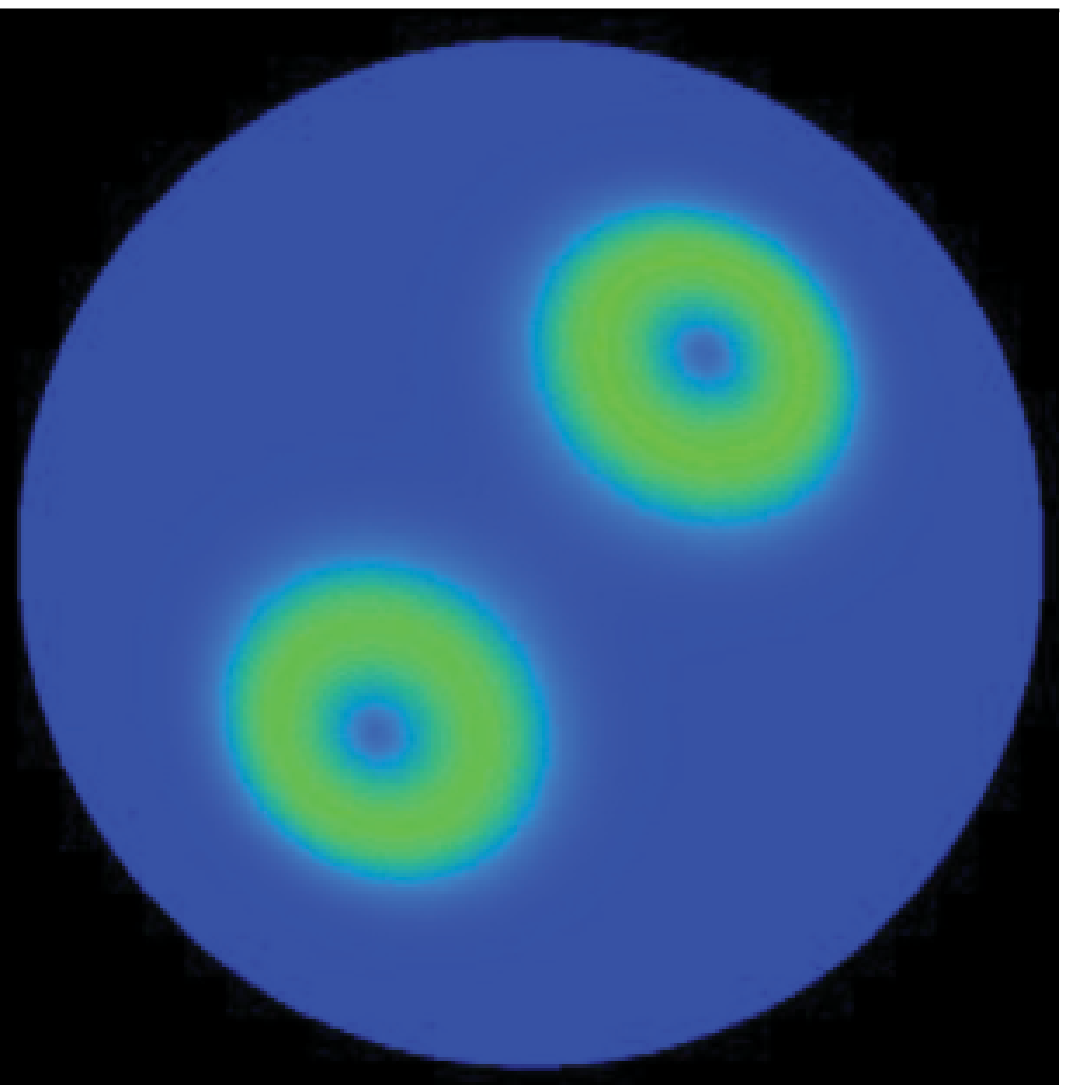}
\includegraphics[scale=0.335]{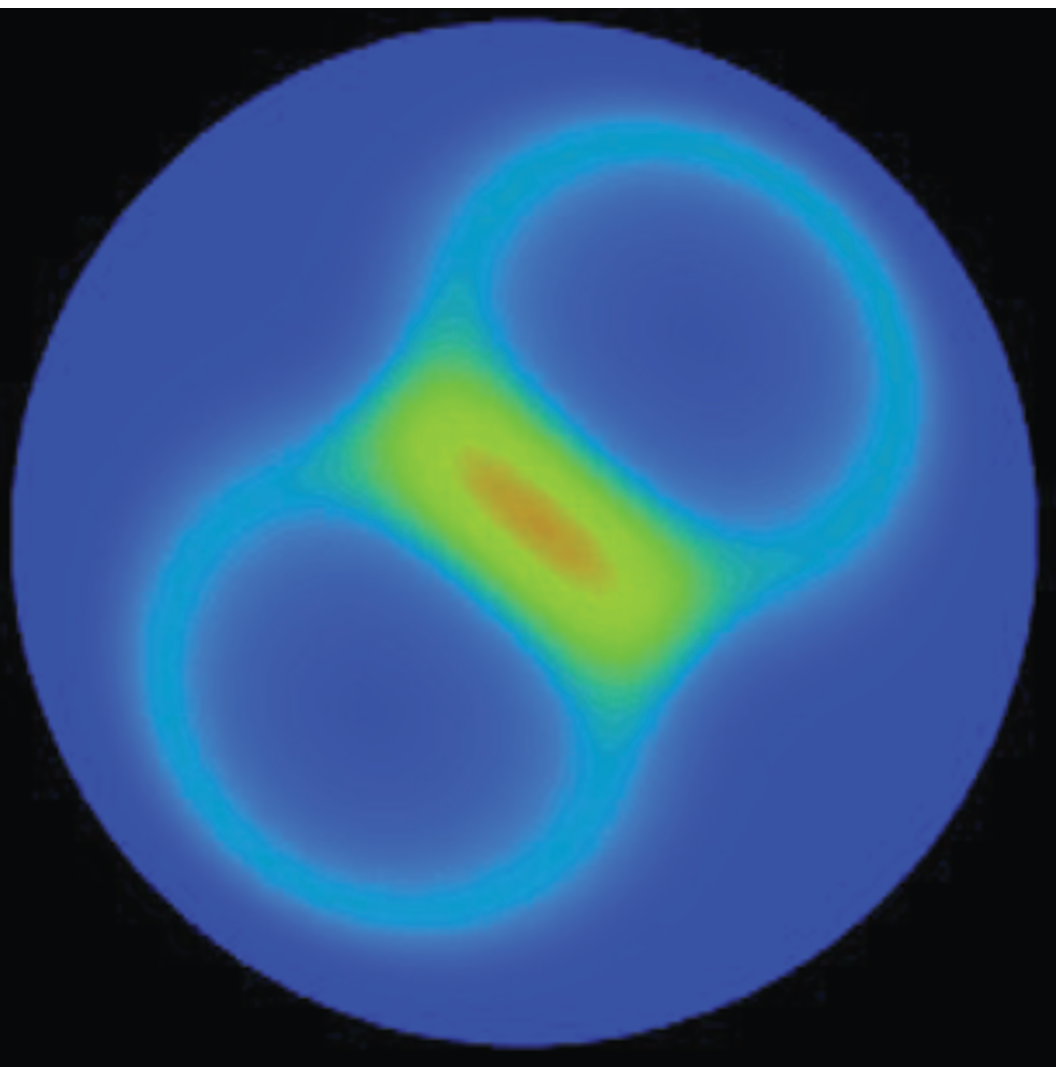}
\includegraphics[scale=0.335]{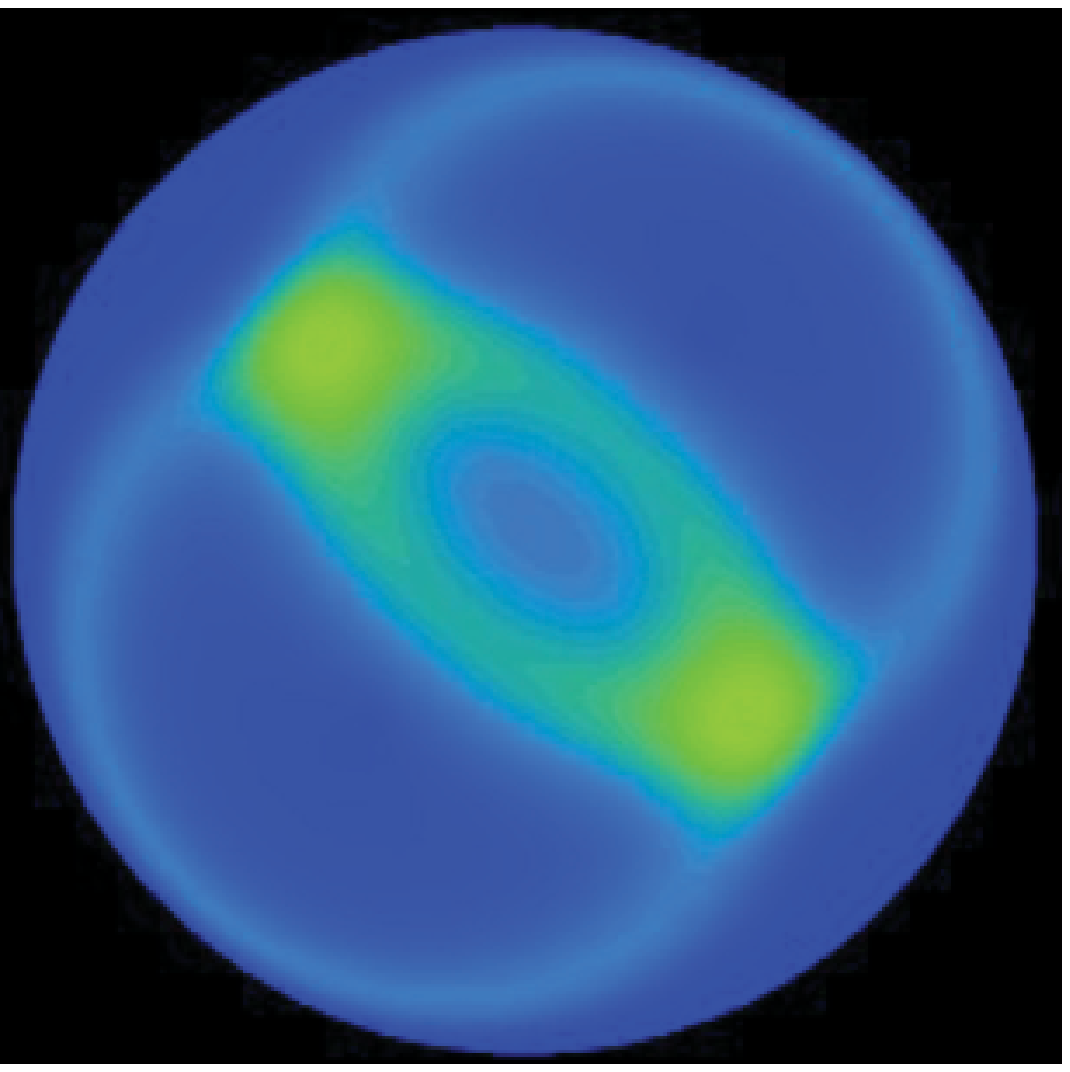}
\end{center}
\begin{center}
{Figure 5:  Simulation of Gaussian pluse interference on sphere by
DEC}
\end{center}
\end{minipage}
$$


\end{document}